\documentclass[12pt]{article}
\usepackage{amsthm,amsfonts,amssymb,amscd}

\begin{document}


\begin{center} \large \bf Affine complements in the projective
space \\ with the trivial group of automorphisms
\end{center}\vspace{0.5cm}

\centerline{A.V.Pukhlikov}\vspace{0.5cm}

\parshape=1
3cm 11cm \noindent {\small \quad\quad\quad \quad\quad\quad\quad
\quad\quad\quad {\bf }\newline We construct a new class of affine
complements ${\mathbb P}^M\setminus S$ with the trivial group of
automorphisms, where $S\subset {\mathbb P}^M$ is a rational
hyper\-surface, $M$ is odd and $M\geqslant 5$.

Bibliography: 19 items.} \vspace{1cm}

\vspace{1cm}

\vspace{1cm}

\section*{Introduction}

{\bf 0.1. Automorphisms of affine complements.} Let $X$ be a
non-singular projective rationally connected variety, $Y\subset X$
a prime ample divisor. There is a natural problem of computing the
group of automorphisms $\mathop{\rm Aut} (X\setminus Y)$ of the
affine variety $X\setminus Y$. This problem can be seen both as a
problem of affine algebraic geometry and as a problem of
birational geometry: the automorphisms of the complement
$X\setminus Y$ are birational self-maps of the variety $X$,
regular on the open set $X\setminus Y$ and transforming it into
itself. Of course, the group of biregular automorphisms
$(\mathop{\rm Aut} X)_Y$, transforming $Y$ into itself, is a
subgroup of the group $\mathop{\rm Aut} (X\setminus Y)$, and we
could say that the problem of computing the latter group is
(practically) solved, if the equality
$$
\mathop{\rm Aut} (X\setminus Y)=(\mathop{\rm Aut} X)_Y
$$
holds, which is of course true if the variety $Y$ is not
birationally ruled (see \cite[Proposition 1]{Pukh2018c}). Since in
any case
$$
\mathop{\rm Aut} (X\setminus Y)\subset \mathop{\rm Bir} X,
$$
the problem of computing the group of automorphisms of an affine
complement is non-trivial only if the group $\mathop{\rm Bir} X$
of birational self-maps of the variety $X$ is sufficiently large,
in particular, infinite. By what was said, it makes sense to
consider only pairs $X\supset Y$, such that $X$ has a large (even
very large) group of birational self-maps, and the ample divisor
$Y$ is a birationally ruled variety. In the present paper, as in
\cite{Pukh2018c}, we consider the case when $X={\mathbb P}^M$ is
the complex projective space and $Y=S$ is a rational hypersurface.
In order to state the main result, let us describe the family of
these hypersurfaces.\vspace{0.3cm}


{\bf 0.2. Rational hypersurfaces.} Fix the integers $e\geqslant 2$
and $m\geqslant e+1$. Set $M=2e+1$. Consider the complex
projective space ${\mathbb P}={\mathbb P}^M$ of dimension
$M\geqslant 5$ with the homogeneous coordinates
$$
(x_0:\dots:x_e:y_0:\dots:y_e)
$$
(the meaning of these non-standard notations will become clear
below). By the symbol $q_{\alpha,\beta}(x_*,y_*)$ we denote a
bi-homogeneous polynomial in these coordinates of degree
$\alpha\in{\mathbb Z}_+$ in $x_*$ and $\beta\in{\mathbb Z}_+$ in
$y_*$. Obviously, such polynomials form a linear subspace of the
space of homogeneous polynomials of degree $\alpha+\beta$ on
${\mathbb P}^M$. The symbol ${\cal P}(m)$ stands for the linear
space of polynomials of the form
$$
q_{m,m+1}(x_*,y_*)+q_{m+1,m}(x_*,y_*).
$$
By the symbol ${\cal P}_{\rm sym}(m)$ we denote the linear
subspace in ${\cal P}(m)$, consisting of polynomials, symmetric
with respect to the replacement $x_i\leftrightarrow y_i$,
$i=0,\dots, e$, that is,
$$
q_{m+1,m}(y_*,x_*)=q_{m,m+1}(x_*,y_*).
$$
Finally, let
$$
P_1=\{y_0=\dots=y_e=0\},\quad P_2=\{x_0=\dots=x_e=0\}
$$
be subspaces in ${\mathbb P}^M$. Obviously, $(x_0:\dots:x_e)$ and
$(y_0:\dots:y_e)$ form natural systems of homogeneous coordinates
on $P_1\cong{\mathbb P}^e$ and $P_2\cong{\mathbb P}^e$. For
$f\in{\cal P}(m)$ we have $P_i\subset S(f)=\{f=0\}$, $i=1,2$, and
moreover, if both components $q_{m,m+1}$ and $q_{m+1,m}$ are not
identically zero, then for $i=1,2$ the equality
$$
\mathop{\rm mult}\nolimits_{P_i}S(f)=m
$$
holds. Let $p_i\in P_i$ be points of general position and
$[p_1,p_2]\subset{\mathbb P}$ the line going through these points.
The restriction of the polynomial $f$ onto that line gives the
divisor
$$
\left(f|_{[p_1,p_2]}\right)=mp_1+mp_2+R(p_1,p_2),
$$
where $R(p_1,p_2)\in S(f)$ is some point. We get a map
$$
R\colon P_1\times P_2\dashrightarrow S(f),
$$
which is birational (the inverse map $R^{-1}\colon
S(f)\dashrightarrow P_1\times P_2$ is $\pi_2\times\pi_1|_{S(f)}$,
where $\pi_i\colon {\mathbb P} \dashrightarrow P_j$ is the linear
projection from the subspace $P_i$, $\{i,j\}=\{1,2\}$). Therefore,
the hypersurface $S(f)$ is rational.\vspace{0.3cm}


{\bf 0.3. Statement of the main result.} Let $\sigma\colon{\mathbb
P}^+\to{\mathbb P}$ be the blow up of the disjoint union $P_1\cup
P_2$ and $\Delta_i=\sigma^{-1}(P_i)$ the non-singular exceptional
divisors,
$$
\Delta_i\cong P_i\times{\mathbb P}^e\cong{\mathbb
P}^e\times{\mathbb P}^e,
$$
where
$$
((x_0:\dots:x_e), (y_0:\dots:y_e))
$$
are bi-homogeneous coordinates on $\Delta_1$ and
$$
((y_0:\dots:y_e), (x_0:\dots:x_e))
$$
are bi-homogeneous coordinates on $\Delta_2$. Here the first group
of coordinates corresponds to $P_i$ and the second one to the
direct factor. Using the same symbols as bi-homogeneous
coordinates on $\Delta_i$, we somewhat abuse our notations,
however, the second group of coordinates corresponds naturally to
precisely the group of coordinates on ${\mathbb P}$ that
complements the coordinates on $P_i$, so that these notations
cannot lead to any misunderstanding. In terms of these systems of
bi-homogeneous coordinates the intersection $S^+\cap\Delta_1$ is
given by the equation $q_{m+1,m}(x_*,y_*)=0$, and the intersection
$S^+\cap\Delta_2$ by the equation $q_{m,m+1}(x_*,y_*)=0$ (where
$S^+\subset{\mathbb P}^+$ is the strict transform of the
hypersurface $S=S(f)$). It is easy to show (see Subsection 1.1),
that for a general polynomial $f\in{\cal P}(m)$ (or ${\cal P}_{\rm
sym}(m)$) the hypersurface $S^+$ and its intersections
$S^+\cap\Delta_i$, $i=1,2$, are non-singular. Denote by the
symbols ${\cal S}$ and ${\cal S}_{\rm sym}$ the linear systems
$$
\{S(f)\,|\, f\in{\cal P}(m)\}\quad\mbox{and} \quad\{S(f)\,|\,
f\in{\cal P}_{\rm sym}(m)\}
$$
of divisors on ${\mathbb P}$, respectively.

Set
$$
{\cal X}\subset{\mathbb P}(H^0({\mathbb P}^e,{\cal O}_{{\mathbb P}^e}(m)))
$$
to be the open $\mathop{\rm Aut}{\mathbb P}^e$-invariant set of
non-singular hypersurfaces of degree $m$ in ${\mathbb P}^e$. Since
the group of automorphisms of a smooth hypersurface of degree $m$
in ${\mathbb P}^e$ is finite \cite{MatMon}, the action of the
$(e+1)^2-1=(e^2+2e)$-dimensional group $\mathop{\rm Aut}{\mathbb
P}^e$ fibres ${\cal X}$ into $(e^2+2e)$-dimensional orbits.

Let $\sigma_S=\sigma|_{S^+}\colon S^+\to S$ be the restriction of
the blow up $\sigma$ onto $S^+$, that is, the blow up of the
singular set $P_1\cup P_2$ on $S$, so that
$$
\sigma^{-1}_S(P_i)=S^+\cap\Delta_i
$$
are the exceptional divisors of that blow up. Assuming that they
are non-singular, we get that for a point $p\in P_i$ of general
position $\sigma^{-1}_S(p)=\sigma^{-1}_S(p)\cap S^+$ is a
non-singular hypersurface of degree $m$ in ${\mathbb P}^e$, that
is, a point of the open set ${\cal X}$, which defines the maps
$$
\xi_i\colon P_i\dashrightarrow{\cal X},\quad i=1,2.
$$
The closure of the image of $P_i$ in ${\cal X}$ we denote by the
symbol $\xi_i(P_i)$. Since $\mathop{\rm dim}P_i=e$ and
$\mathop{\rm dim}{\cal X}={m+e\choose e}-1\geqslant{2e+1\choose
e}-1$, the following condition of general position for the
polynomial $f$ (or the hypersurface $S$) makes sense:

$(*)$ for some point $p\in P_i$, for which the hypersurface
$\sigma^{-1}_S(p)\subset{\mathbb P}^e$ is non-singular, the image
$\xi_i(P_i)$ intersects the $\mathop{\rm Aut}{\mathbb P}^e$-orbit
of the point $\xi_i(p)$ at this point only, $i=1,2$.

Obviously, if the hypersurface $S$ satisfies the condition $(*)$,
then this condition holds for every point of general position
$p\in P_i$.

Now we can state the main result of the present paper.

{\bf Theorem 0.1.}  {\it Assume that the strict transform $S^+$ of
the hypersurface $S\in{\cal S}$ is non-singular and moreover, the
exceptional divisors $\sigma^{-1}_S(P_i)$ are non-singular,
$i=1,2$.

{\rm (i)} If $S'\subset{\mathbb P}$ is an irreducible hypersurface
of degree $2m+1$, such that there is an isomorphism
$$
\chi\colon{\mathbb P}\backslash S\to{\mathbb P}\backslash S'
$$
of affine algebraic varieties ${\mathbb P}\backslash S$ and
${\mathbb P}\backslash S'$, then $\chi$ is the restriction onto
${\mathbb P}\backslash S$ of a projective automorphism
$\chi_{\mathbb P}\in\mathop{\rm Aut}{\mathbb P}$, transforming $S$
into $S'$, so that $S$ and $S'$ are isomorphic.}

{\rm (ii)} {\it If $S$ satisfies the condition $(*)$ and the
closed subsets $\xi_1(P_1)$ and $\xi_2(P_2)$ in ${\cal X}$ are not
equal, then the group of automorphisms $\mathop{\rm Aut}({\mathbb
P}\backslash S$) is trivial.}

(iii) {\it If $S\in{\cal S}_{\rm sym}$ satisfies the condition
$(*)$, then
$$
\mathop{\rm Aut}({\mathbb P}\backslash
S)=\langle\tau\rangle\cong{\mathbb Z}/2{\mathbb Z},
$$
where the projective involution $\tau$ is defined by the
equalities}
$$
\tau^*x_j=y_j,\quad \tau^*y_j=x_j,\quad
j=0,1,\dots,e.
$$\vspace{0.3cm}


{\bf 0.4. The structure of the paper.} The paper is organized in
the following way. In \S 1 we prepare to the proof of the main
theorem. In Subsection 1.1 we consider the singularities of the
hypersurface $S$, justifying the assumptions of the theorem. In
Subsection 1.2 the proof starts: in the assumption that the
isomorphism of the affine varieties ${\mathbb P}\setminus S$ and
${\mathbb P}\setminus S'$ is the restriction of a birational, but
not a biregular, automorphism of the projective space, we resolve
the singularities of the latter map and define the key objects and
parameters, characterizing this map, and prove an analog of
\cite[Proposition 2]{Pukh2018c} that gives a relation between
these parameters. Besides, we note that the result of
\cite{Pukh2018c} can be made somewhat stronger (no changes in the
proof are needed).

In \S 2 we prove the main theorem. It is done in a few steps. To
begin with, we consider the sequence of blow ups, resolving the
exceptional divisor $T$, corresponding to the hypersurface $S'$
(Subsection 2.1). After that, we exclude one by one all options
for the centre of the divisor $T$ on ${\mathbb P}$: when it is not
contained in the singular set of the hypersurface $S$ and when it
is contained in one of the subspaces $P_1$ or $P_2$ (Subsections
2.2--2.4). This proves the claim (i) of the main theorem. In
Subsection 2.5 we prove the claims (ii) and (iii), which completes
the proof of Theorem 0.1.\vspace{0.1cm}


{\bf 0.5. Remarks.} Studying the groups of automorphisms of affine
algebraic varieties is a classical and established topic in
algebraic geometry. As usual, in the dimension 2 quite a lot is
known: \cite{GizDan1975,GizDan1977} are dealing with groups of
automorphisms of affine surfaces (as far as the author knows,
V.A.Iskovskikh also had a proof of the theorem on the group of
automorphisms of the affine plane ${\mathbb A}^2$, which, it
seems, has never been published). The two-dimensional case is
still a field of active research, see, for instance,
\cite{DubLamy2015,FurLamy2010,CanLamy2006,Lamy2005}.

In the dimension 3 and higher we know a lot less. There are just a
few papers with complete results on the groups of automorphisms of
affine varieties (in the cases when the problem is not trivial),
see \cite{Pukh2018c,CheltsovDuboulozPark2018}. Such classical
problems as computing the group of automorphisms of the complement
in ${\mathbb P}^3$ to a non-singular cubic surface (see
\cite{Giz2005}) or the group of automorphisms of the affine space
${\mathbb A}^3$ remain open.

The problems of describing the groups of automorphisms of affine
varieties are on the border with such areas of algebraic geometry
as describing the groups of birational self-maps and their
subgroups and the theory of birational rigidity. A considerable
recent progress in these areas (see, for instance,
\cite{PrShr2016,PrShr2017,PrShr2018a,PrShr2018b,PrzShr2021,
Pukh2022b,Prokh2023}), gives a basis for some optimism in respect
of the groups of automorphisms of affine varieties, too.


\section{Preliminary constructions}

In this section we consider the singularities of a general
hypersurface $S$ in more detail (Subsection 1.1) and start the
proof of the main theorem (Subsection 1.2).\vspace{0.3cm}

{\bf 1.1. Singularities of the hypersurface $S$.} Before starting
to prove Theorem 0.1, let us consider in more detail the linear
systems ${\cal S}$ and ${\cal S}_{\rm sym}$, defined in Subsection
0.3. For their strict transforms on ${\mathbb P}^+$ we have the
equalities
$$
{\cal S}^+=\{\sigma^*S(f)-m\Delta_1-m\Delta_2\,|\, f\in{\cal
P}(m)\}
$$
and
$$
{\cal S}^+_{\rm sym}=\{\sigma^*S(f)-m\Delta_1-m\Delta_2\,|\,
f\in{\cal P}_{\rm sym}(m)\},
$$
respectively.

{\bf Proposition 1.1.} {\it} The following equality holds:
$\mathop{\rm Bs}{\cal S}^+=\mathop{\rm Bs}{\cal S}^+_{\rm
sym}=\emptyset$.

{\bf Proof.} Since the polynomials $q_{m+1,m}$ and $q_{m,m+1}$ are
arbitrary, these linear systems have no base points on
$\Delta_1\cup\Delta_2$, whereas for any point $p\in{\mathbb
P}\backslash(P_1 \cup P_2)$ and a general polynomial $f$ we have
$f(p)\neq 0$. Q.E.D. for the proposition.

It follows that
$$
\mathop{\rm Bs}{\cal S}=\mathop{\rm Bs}{\cal S}_{\rm sym}=P_1\cup
P_2
$$
and for a general divisor $S\in{\cal S}$ (or ${\cal S}_{\rm sym}$)
its strict transform $S^+\subset{\mathbb P}^+$ is non-singular. It
is also obvious that for a general divisor $S\in{\cal S}$ (or
${\cal S}_{\rm sym}$) the intersections $S^+\cap\Delta_i$ are
non-singular divisors of bi-degree $(m+1,m)$ on
$\Delta_i\cong{\mathbb P}^e\times{\mathbb P}^e$. Since $m\geqslant
e+1$, they are not rationally connected (and even uniruled)
varieties.\vspace{0.3cm}


{\bf 1.2. The birational self-map $\chi_{\mathbb P}$.} Let us
start the proof of Theorem 0.1. Let
$$
\chi\colon{\mathbb P}\backslash S\to{\mathbb P}\backslash S'
$$
be a non-trivial isomorphism of affine varieties, that is, the
corresponding birational self-map $\chi_{\mathbb P}\in\mathop{\rm
Bir}{\mathbb P}$ is not a biregular automorphism of the projective
space. Consider the composition
$$
\chi_{\mathbb P}\circ\sigma\colon{\mathbb
P}^+\dashrightarrow{\mathbb P}.
$$
Let $\varphi\colon\widetilde{\mathbb P}\to{\mathbb P}^+$ be a
resolution of singularities of that map (we will see below that
$\chi_{\mathbb P}\circ\sigma$ can not be regular, but for the
moment we assume that if $\chi_{\mathbb P}\circ\sigma$ is regular,
then $\varphi$ is the identity map). Let
$$
T=(\chi_{\mathbb
P}\circ\sigma\circ\varphi)^{-1}_*S'\subset\widetilde{\mathbb P}
$$
be the strict transform of the hypersurface $S'$. Since
$\chi_{\mathbb P}\not\in\mathop{\rm Aut}{\mathbb P}$, we see that
the image $(\sigma\circ\varphi)(T)\subset{\mathbb P}$ is not the
hypersurface $S$, so that $T$ is a
$(\sigma\circ\varphi)$-exceptional divisor. Set
$$
B_{\mathbb P}=(\sigma\circ\varphi)(T) \subset{\mathbb P}.
$$
This is a subvariety of codimension $\geqslant 2$. Set also
$B=\varphi(T)\subset{\mathbb P}^+$. It is possible that
$B=\Delta_1$ or $\Delta_2$, but if this is not the case, then
$B\subset{\mathbb P}^+$ is a subvariety of codimension $\geqslant
2$. Obviously, $B_{\mathbb P}\subset S$ and
$$
B\subset S^+\cup\Delta_1\cup\Delta_2,
$$
where $\sigma(B)=B_{\mathbb P}$. Let
$$
a=a(T,{\mathbb P})\geqslant 1\quad\mbox{and}\quad a_+=a(T,{\mathbb
P}^+)\geqslant 0
$$
be the discrepancies and
$$
b=\mathop{\rm ord}\nolimits_T\varphi^*\sigma^*S\quad
\mbox{and}\quad b_+=\mathop{\rm ord}\nolimits_T\varphi^*S^+.
$$
If $T=\Delta_1$ or $\Delta_2$, then $a_+=b_+=0$, otherwise
$a_+\geqslant 1$. If $B\not\subset\Delta_1\cup\Delta_2$, then
$a=a_+$ and $b=b_+$; in any case $b\geqslant 1$. If
$B\subset\Delta_1\cup\Delta_2$, then since
$\Delta_1\cap\Delta_2=\emptyset$, there is precisely one value
$i\in\{1,2\}$, such that $B\subset\Delta_i$, and in that case set
$$
\varepsilon=\mathop{\rm ord}\nolimits_T\varphi^*\Delta_i.
$$
If $B\not\subset\Delta_1\cup\Delta_2$, then set $\varepsilon=0$.

We get the obvious equalities
$$
a=a_++\varepsilon e,\quad b=b_++\varepsilon m.
$$

{\bf Proposition 1.2.} {\it The following equality holds:}
$$
(M+1)b=a(2m+1)+(M+1).
$$

{\bf Proof} repeats the proof of the second equality (3.1) in
\cite{Pukh2018c} almost word for word. Namely, let
$\widetilde{S}\subset\widetilde{\mathbb P}$ be the strict
transform of the hypersurface $S$ and
$\widetilde{K}=K_{\widetilde{\mathbb P}}$ the canonical class. We
get the equalities
$$
\widetilde{K}=-(M+1)H+aT+(\dots),
$$
$$
\widetilde{S}\sim(2m+1)H-bT+(\dots),
$$
where for simplicity we write $H$ instead of $\varphi^*\sigma^*H$
and in the brackets we have integral linear combinations of
$(\sigma\circ\varphi)$-exceptional divisors that are different
from $T$. Their images with respect to the morphism $\chi_{\mathbb
P}\circ\sigma\circ\varphi$ are of codimension $\geqslant 2$ in
${\mathbb P}$ and for that reason a general line $L\subset{\mathbb
P}$ does not meet them. Furthermore, $\chi_{\mathbb
P}\circ\sigma\circ\varphi(\widetilde{S})$ also is of codimension
$\geqslant 2$ in ${\mathbb P}$. Therefore, for the strict
transform $\widetilde{L}\subset\widetilde{\mathbb P}$ the
following equalities hold:
$$
(\widetilde{L}\cdot\widetilde{S})=0,\quad
(\widetilde{L}\cdot\widetilde{K})=-(M+1),
$$
and $(\widetilde{L}\cdot T)=(L\cdot S')=2m+1$. Set
$d=(\widetilde{L}\cdot H)$, it is the degree of the curve
$\sigma\circ\varphi(\widetilde{L})\subset{\mathbb P}$ in the usual
sense. As a result we get
$$
(2m+1)d-b(2m+1)=0,
$$
so that $d=b$, and
$$
-(M+1)d+a(2m+1)=-(M+1),
$$
which completes the proof of the proposition.

Now let us show that
$T\neq\widetilde{\Delta_1},\widetilde{\Delta_2}$, so that the map
$\chi_{\mathbb P}\circ\sigma$ can not be regular. Assume that
$T=\widetilde{\Delta_1}$. In the notations of the proof of
Proposition 1.2 we get $d=b=m$, so that the curve
$\sigma\circ\varphi(\widetilde{L})$ is a rational curve of degree
$m$ in ${\mathbb P}$, the strict transform of which
$\varphi(\widetilde{L})\subset{\mathbb P}^+$ meets $\Delta_1$ at
$(2m+1)$ points of general position. This is impossible: take a
general hyperplane $\Xi\supset P_1$, then
$(\Xi\cdot[\sigma\circ\varphi(\widetilde{L})])=m$ and
$(\Xi^+\cdot(\varphi(\widetilde{L}))=-m-1$, since
$\sigma^*\Xi=\Xi^++\Delta_1$, where $\Xi^+\subset{\mathbb P}^+$ is
the strict transform.

{\bf Remark 1.1.} In \cite{Pukh2018c} for the hypersurface $S$ a
hypersurface of degree $\mathop{\rm deg}S\geqslant M+1$ with a
unique singular point $o\in S$ of multiplicity $\mathop{\rm
deg}S-1$, resolved by one blow up, was considered. The notations
of \cite{Pukh2018c} are close to the notations of the present
paper, and where they are identical, they can be used in both
contexts. In \cite{Pukh2018c} in the proof of regularity of the
map $\chi_{\mathbb P}$ we missed the case, similar to the one that
we have just considered, when the hypersurface $S'$ is transformed
to the the exceptional divisor of the blow up of the point $o$ on
${\mathbb P}={\mathbb P}^M$. That case is easy to exclude by the
arguments similar to our arguments above: the strict transform of
a general line with respect to $\chi_{\mathbb P}$ is a rational
curve of degree $\mathop{\rm mult}_oS=\mathop{\rm deg}S-1$ with
the multiplicity $\mathop{\rm deg}S'=\mathop{\rm deg}S$ at the
point $o$, which is impossible. Note that the claim of Theorem 2
in \cite{Pukh2018c} holds for any irreducible hypersurface
$S'\subset{\mathbb P}$ of degree $\mathop{\rm deg}S$, since the
proof does not use any other information about $S'$. The
assumptions about that hypersurface made in \cite[Subsection
1.1]{Pukh2018c} are unnecessary.


\section{Proof of the main result}

In this section we prove Theorem 0.1. To begin with, we consider
the sequence of blow ups, resolving the exceptional divisor $T$,
and certain invariants, related to this sequence (Subsection 2.1).
Then we exclude one by one the non-singular case, when
$\sigma(B)\not\subset \mathop{\rm Sing} S$ (Subsection 2.2) and
the singular case, when $\sigma(B)\subset P_1\cup P_2$
(Subsections 2.3 and 2.4), which completes the proof of the part
(i) of Theorem 0.1. In Subsection 2.5 we prove the claims (ii) and
(iii).\vspace{0.3cm}

{\bf 2.1. The resolution of the exceptional divisor $T$.} As we
saw in \S 1, $T$ is a divisor over ${\mathbb P}^+$, since
$T\neq\Delta_1,\Delta_2,S^+$. Its centre on ${\mathbb P}^+$ is an
irreducible subvariety $\varphi(T)\subset
S^+\cup\Delta_1\cup\Delta_2$ of codimension $\geqslant 2$. Let
$$
\varphi_{i,i-1}\colon X_i\to X_{i-1},
$$
$i=1,\dots,N$, be the resolution of the divisor $T$ in the sense
of \cite[Chapter 2]{Pukh13a}, that is, $X_0={\mathbb P}^+$, the
birational morphism $\varphi_{i,i-1}$ is the blow up of the centre
of $T$ on $X_{i-1}$, which we denote by the symbol $B_{i-1}$,
where $\mathop{\rm codim}(B_{i-1}\subset X_{i-1})\geqslant 2$ for
$i=1,\dots,N$, with the exceptional divisor
$E_i=\varphi^{-1}_{i,i-1}(B_{i-1})$ and $E_N\subset X_N$ is a
realization of the discrete valuation, corresponding to $T$, that
is, the birational map
$$
\varphi^{-1}_{N,N-1}\circ\dots\circ\varphi^{-1}_{1,0}\circ\varphi
\colon\widetilde{{\mathbb P}}\dashrightarrow X_N
$$
is regular at the general point of $T$ and transforms $T$ into
$E_N$. We use the standard notations: for $j>i$ the composition
$$
\varphi_{i+1,i}\circ\dots\circ\varphi_{j,j-1}\colon X_j\to X_i
$$
is denoted by the symbol $\varphi_{j,i}$, and the strict transform
of a subvariety $R\subset X_i$ on $X_j$, if it is well defined, is
denoted by the symbol $R^j$. Note that the blown up subvarieties
$B_{i-1}$, generally speaking, can be singular, however
$$
B_{i-1}\not\subset\mathop{\rm Sing}X_{i-1}
$$
for all $i=1,\dots,N$, and over a general point of $B_{i-1}$ the
variety $X_i$ and the exceptional divisor $E_i$ are non-singular.
By construction, for $j>i$ we have
$$
\varphi_{j-1,i-1}(B_{j-1})=B_{i-1}.
$$
On the set of exceptional divisors $E_1,\dots,E_N$ there is a
natural structure of an oriented graph: an oriented edge (an
arrow) goes from $E_i$ to $E_j$ (which is written down as $i\to
j$), if and only if $i>j$ and $B_{i-1}\subset E^{i-1}_j$. If $E_i$
and $E_j$ are not joined by an oriented edge, we write
$i\nrightarrow j$.

The number of paths from the vertex $E_i$ to $E_j$ for $i\neq j$
in the oriented graph described above we denote by the symbol
$p_{ij}$. Obviously, for $i\neq j$ we have $p_{ij}\geqslant 1$ if
and only if $i>j$, and for $i<j$ we have $p_{ij}=0$. For
convenience we set $p_{ii}=1$. Obviously, for $i>j$ we have the
equality
$$
p_{ij}=\mathop{\rm ord}\nolimits_{E_i}\varphi^*_{i,j}E_j.
$$
(It justifies the convention $p_{ii}=1$.) Besides, in order to
simplify the notations, we set $p_i=p_{Ni}$. (See \cite[Chapter
2]{Pukh13a} or \cite[Chapter 2]{Pukh07a} for the details.)
Finally, for $i=1,\dots,N$ set
$$
\delta_i=\mathop{\rm codim}(B_{i-1}\subset X_{i-1})-1.
$$
Now we get the formula for the discrepancy
$$
a_+=\sum^N_{i=1}p_i\delta_i.
$$
Set also
$$
k=\mathop{\rm max}\{i\,|\, 1\leqslant i\leqslant N, B_{i-1}\subset
S^{i-1}\},
$$
where for simplicity of notations we write $S^{i-1}$ instead of
$(S^+)^{i-1}$; if the set in the right hand side is empty, we set
$k=0$. Since the hypersurface $S^+$ is non-singular, for
$i=1,\dots,k$ the hypersurface $S^{i-1}$ is non-singular at the
general point of $B_{i-1}$, so that we get the formula
$$
b_+=\sum^k_{i=1}p_i,
$$
which remains true for $k=0$, since in that case
$b_+=0$.\vspace{0.3cm}


{\bf 2.2. Exclusion of the non-singular case.} We say that the
non-singular case takes place, if
$$
B_{\mathbb P}=\sigma(B)\not\subset\mathop{\rm Sing}S=P_1\cup P_2,
$$
that is, $B\not\subset\Delta_1\cup\Delta_2$.

{\bf Proposition 2.1.} {\it The non-singular case does not
realize: $B\subset\Delta_i$ for some} $i\in\{1,2\}$.

{\bf Proof.} Assume the converse: $B_{\mathbb
P}\not\subset\mathop{\rm Sing}S$. Then $a=a_+$ and $b=b_+$, so
that from Proposition 1.2, taking into account the explicit
formulas above, we get the equality
$$
(M+1)\sum^k_{i=1}p_i=(2m+1)\sum^N_{i=1}p_i\delta_i+(M+1).
$$
Since $\delta_i\geqslant 1$, $2m+1\geqslant 2e+3\geqslant M+2$ and
$k\leqslant N$, this equality is impossible: the right hand side
is strictly higher than the left hand one. This contradiction
proves the proposition. Q.E.D.

Therefore, $B\subset\Delta_1\cup\Delta_2$. Since
$\Delta_1\cap\Delta_2=\emptyset$, the subvariety $B$ is contained
in precisely one of the divisors $\Delta_1,\Delta_2$, which we now
denote by the symbol $\Delta$. Set
$$
l=\mathop{\rm max}\{i\,|\, 1\leqslant i\leqslant N,
B_{i-1}\subset\Delta^{i-1}\}.
$$
By construction, $l\geqslant 1$. Since the subvariety $\Delta$ is
non-singular, we get the formula
$$
\varepsilon=\mathop{\rm
ord}\nolimits_T\varphi^*\Delta=\sum^l_{i=1}p_i,
$$
so that Proposition 1.2 gives the equality
\begin{equation}\label{24.02.24.1}
(M+1)\left[\sum^k_{i=1}p_i+m\sum^l_{i=1}p_i\right]=(2m+1)
\left[\sum^N_{i=1}p_i\delta_i+e\sum^l_{i=1}p_i\right]+(M+1).
\end{equation}\vspace{0.3cm}


{\bf 2.3. Exclusion of the singular case for $l\leqslant k$.} The
following fact is true (here we do not assume that $l\leqslant
k$).

{\bf Proposition 2.2.} {\it Assume that $k\geqslant 1$. Then for
any $i\leqslant\mathop{\rm min}\{l,k\}$, $i\geqslant 1$, the
inequality} $\delta_i\geqslant 2$ {\it holds.}

{\bf Proof.} By construction, the subvariety $B_{i-1}$ is
contained in $E_{i-1}$, and by the inequalities $i\leqslant l$ and
$i\leqslant k$ it is also contained in the strict transforms
$\Delta^{i-1}$ and $S^{i-1}$. By the assumption of Theorem 0.1 the
hypersurfaces $S^+$ and $\Delta$ meet everywhere transversally,
which implies that the divisor $E_{i-1}+\Delta^{i-1}+S^{i-1}$ has
normal crossings at a general point of the subvariety $B_{i-1}$,
that is,
$$
E_{i-1}\cap\Delta^{i-1}\cap S^{i-1}
$$
is a non-singular subvariety of codimension 3 in $X_{i-1}$. For
that reason, the more so we get $\mathop{\rm codim}(B_{i-1}\subset
X_{i-1})\geqslant 3$, which completes the proof of the
proposition.

Assume now that $l\leqslant k$. In that case the equality
(\ref{24.02.24.1}) can be re-written as the equality, in the left
hand side of which we have the sum
$$
\sum^l_{i=1} ((M+1)(m+1)-(2m+1)(\delta_i+e))p_i,
$$
and in the right hand side we have the expression
$$
\sum^k_{i=l+1}((2m+1)\delta_i-M-1)
p_i+(2m+1)\sum^N_{i=k+1}p_i\delta_i+(M+1).
$$
Since $M=2e+1$, $m\geqslant e+1$ and $\delta_i\geqslant 1$, the
right hand side is strictly positive. However, for $i\leqslant l$
by Proposition 2.2 we have $\delta_i\geqslant 2$ and it is easy to
check that
$$
(M+1)(m+1)-(2m+1)(e+2)=e-2m<0.
$$
Thus the left hand side is strictly negative. We obtained a
contradiction that excludes the singular case for $l\leqslant
k$.\vspace{0.3cm}


{\bf 2.4. Exclusion of the singular case for $l>k$.} In that case
the equality (\ref{24.02.24.1}) can be re-written as the equality,
in the left hand side of which we have the sum
$$
\sum^k_{i=1}((M+1)(m+1)-(2m+1)(\delta_i+e))p_i
$$
(if $k=0$, then it is equal to zero), and in the right hand side
the expression
$$
\sum^l_{i=k+1}((2m+1)(\delta_i+e)-(M+1)m)p_i+
(2m+1)\sum^N_{i=l+1}p_i\delta_i+(M+1).
$$
Thus the left hand side is non-positive (strictly negative for
$k\geqslant 1$, as we saw in Subsection 2.3). Since
$\delta_i\geqslant 1$ and the equality
$$
(2m+1)(e+1)-(2e+2)m=e+1
$$
holds, we get a contradiction again, excluding the singular case
for $l>k$. This completes the proof of the claim (i) of Theorem
0.1.\vspace{0.3cm}


{\bf 2.5. Projective automorphisms, preserving the hypersurface
$S$.} Let us prove the claim (ii) of Theorem 0.1. We know already
that $\chi_{\mathbb P}\in\mathop{\rm Aut}{\mathbb P}$ is a
projective automorphism, preserving $S$. Since $\mathop{\rm
Sing}S=P_1\cup P_2$, we have that $\chi_{\mathbb P}(P_1)$ is
either $P_1$ or $P_2$. However, by assumption
$\xi_1(P_1)\neq\xi_2(P_2)$, so that $\chi_{\mathbb P}(P_1)=P_1$
and $\chi_{\mathbb P}(P_2)=P_2$. By the condition $(*)$ for a
point of general position $p\in P_i$ the fibre of the exceptional
divisor $\Delta_i$ over $p$ is not isomorphic to the fibre of
$\Delta_i$ over any other point of the subspace $P_i$. Therefore
$$
\chi_{\mathbb P}|_{P_i}=\mathop{\rm id}\nolimits_{P_i}.
$$
Now let us consider points of general position $p_1\in P_1$ and
$p_2\in P_2$ and the $[p_1,p_2]$ passing through these points. The
automorphism $\chi_{\mathbb P}$ preserves that line and its
restriction $\chi_{\mathbb P}|_{[p_1,p_2]}$ preserves the three
points: $p_1,p_2$ and the points
$$
([p_1,p_2]\cap S)\backslash\{p_1,p_2\}.
$$
Therefore, $\chi_{\mathbb P}|_{[p_1,p_2]}=\mathop{\rm
id}_{[p_1,p_2]}$, so that $\chi_{\mathbb P}=\mathop{\rm
id}_{\mathbb P}$, as we claimed. Q.E.D. for the part (ii).

Let us show (iii). Here we argue in a similar way with the only
difference: $\chi_{\mathbb P}$ can swap the subspaces $P_1,P_2$.
If $\chi_{\mathbb P}(P_i)=P_i$, then we can repeat the proof above
word for word and conclude that $\chi_{\mathbb P}=\mathop{\rm
id}_{\mathbb P}$. If $\chi_{\mathbb P}(P_1)=P_2$, then, taking
into account that the hypersurface $S$ is invariant with respect
to $\tau$, we consider the composition $\tau\chi_{\mathbb
P}\in\mathop{\rm Aut}{\mathbb P}$. This transformation preserves
the subspaces $P_1,P_2$, so that $\tau\chi_{\mathbb P}=\mathop{\rm
id}_{\mathbb P}$ as above and $\chi_{\mathbb P}=\tau$. Q.E.D. for
our theorem.


\begin{flushleft}
Department of Mathematical Sciences,\\
The University of Liverpool
\end{flushleft}

\noindent{\it pukh@liverpool.ac.uk}


\begin{thebibliography}{99}

\bibitem{Pukh2018c} Pukhlikov A.~V., Automorphisms of certain
affine complements in projective space. Sbornik: Mathematics {\bf
209} (2018), No. 2, 276--289.

\bibitem{MatMon} Matsumura H. and Monsky P., On the automorphisms
of hypersurfaces. J. Math. Kyoto Univ. {\bf 3} (1964), 347--361.

\bibitem{GizDan1975} Gizatullin M.~Kh. and Danilov V.~I.,
Automorphisms of affine surfaces. I. Izvestiya:
Mathematics {\bf 39} (1975), No. 3, 523--565.

\bibitem{GizDan1977} Gizatullin M.~Kh. and Danilov V.~I.,
Automorphisms of affine surfaces. II. Izvestiya:
Mathematics {\bf 41} (1977), No. 1, 54--103.

\bibitem{DubLamy2015} Dubulouz A. and Lamy S., Automorphisms of
open surfaces with irreducible boundary. Osaka Math. J. {\bf 52}
(2015), No. 3, 747--791.

\bibitem{FurLamy2010} Furter J.-Ph. and Lamy S., Normal subgroup
generated by a plane polynomial automorphism. Transform. Groups
{\bf 15} (2010), No. 3, 577--610.

\bibitem{CanLamy2006} Cantat S. and Lamy S., Groupes
d'automorphismes polynomiaux du plan. Geom. Dedicata {\bf 123}
(2006), 201--221.

\bibitem{Lamy2005} Lamy S., Sur la structure du groupe
d'automorphismes de certaines surfaces affines. Publ. Math. {\bf
49} (2005), No. 1, 3--20.

\bibitem{CheltsovDuboulozPark2018} Cheltsov I., Dubouloz A. and
Park J., Super-rigid affine Fano varieties. Compos. Math. {\bf
154} (2018), No. 11, 2462--2484.

\bibitem{Giz2005} Freudenburg G. and Russell P., Open problems in
affine algebraic geometry, in: Affine Algebraic Geometry, Contemp.
Math. {\bf 369}, Amer. Math. Soc., Providence, RI 2005, 1--30.

\bibitem{PrShr2016} Prokhorov Yu. and Shramov C., Jordan property
for Cremona groups, Amer. J. Math. {\bf 138} (2016), No. 2,
403--418.

\bibitem{PrShr2017} Prokhorov Yu. and Shramov C., Jordan constant
for Cremona group of rank 3, Mosc. Math. J. {\bf 17} (2017), No.
3, 457--509.

\bibitem{PrShr2018a} Prokhorov Yu. and Shramov C., $p$-subgroups
in the space Cremona group, Math. Nachr. {\bf 291} (2018), No.
8-9, 1374--1389.

\bibitem{PrShr2018b} Prokhorov Yu. and Shramov C., Finite groups
of birational self-maps of three-folds, Math. Res. Lett. {\bf 25}
(2018), No. 3, 957-972.

\bibitem{PrzShr2021} Przyjalkowski V. and Shramov C., On
automorphisms of quasi-smooth weighted complete intersections,
Sbornik: Mathematics {\bf 212} (2021), No. 3, 374--388.

\bibitem{Pukh2022b} Pukhlikov A.~v., Effective
results in the theory of birational rigidity, Russian Math.
Surveys {\bf 77} (2022), No. 2, 301--354.

\bibitem{Prokh2023} Prokhorov Yu., Finite groups of birational
transformations, in: European Congress of mathematics, EMS Press,
Berlin, 2023, 413--437.

\bibitem{Pukh13a} Pukhlikov Aleksandr, Birationally Rigid
Varieties. Mathematical Surveys and Monographs {\bf 190}, AMS,
2013.

\bibitem{Pukh07a} Pukhlikov A.~V., Birationally rigid
varieties. I. Fano varieties, Russian Math. Surveys. {\bf 62}
(2007), No. 5, 857--942.

\end{thebibliography}
\end{document}